\documentclass[12pt,twoside,final,reqno]{amsart}
\usepackage{amssymb,amsmath,amscd,amsfonts,amsthm,mathrsfs}
\usepackage{verbatim,paralist}
\usepackage{enumerate}
\RequirePackage{ifthen}
\usepackage[latin1]{inputenc}
\usepackage[T1]{fontenc}
\usepackage[english]{babel}
\usepackage{babel}
 \usepackage{tikz}
 \usepackage{float}
\usepackage{latexsym}
\usepackage{amsmath}
\usepackage{amsfonts}
\usepackage{amssymb}
\usepackage{latexsym}
\usepackage{color}
\usepackage{fancyhdr}

\usepackage[latin1]{inputenc}

\newtheorem{theorem}{Theorem}[section]
\newtheorem{Proposition}[theorem]{Proposition}
\newtheorem{Lemma}[theorem]{Lemma}
\newtheorem{Remark}[theorem]{Remark}
\newtheorem{example}[theorem]{Example}
\newtheorem{Definition}[theorem]{Definition}

\title {Fractional Powers of $3\times3$ Block Operator Matrices: Applications to PDEs}
\author{Hatem Baloudi}
\address{Hatem Baloudi, Department of Mathematics, Faculty of Sciences of Gafsa, University of Gafsa, 2112 Zarroug,
Tunisia}
\author{Mohamed Ali Dbeibia}
\address{Mohamed Ali Dbeibia, Departement of Mathematics Faculty of Sciences of Sfax, University of Sfax, Route de Soukra km 3.5, B.P. 1171, 3000 Sfax, Tunisia.}
\author{Aref JERIBI}
\address{Aref Jeribi, Department of Mathematics, Faculty of Sciences of Sfax, University of Sfax, 3000 Sfax,
Tunisia}

\begin{document}

\maketitle
\begin{abstract}In this paper, we investigate the fractional powers of
$3\times3$ block operator matrices, with a particular focus on their applications to partial differential equations (PDEs). We develop a comprehensive theoretical framework for defining and calculating fractional powers of positive operators and extend these results to block operator matrices. Various methods, including alternative formulas, change of variables, and the second resolvent identity, are employed to obtain explicit expressions for fractional powers. The results are applied to systems of PDEs, demonstrating the relevance and effectiveness of the proposed approaches in modeling and analyzing complex dynamical systems. Examples are provided to illustrate the theoretical findings and their applicability to concrete problems.\end{abstract}

\tableofcontents

\section{Introduction}
In this paper, we undertake a comprehensive exploration of the fractional powers of $3\times3$ block operator matrices and their applications 
to partial differential equations (PDEs). The concept of fractional powers of operators, initially developed in the context of functional analysis, 
 has emerged as a powerful tool in the study of evolution equations, 
  abstract differential equations, and complex dynamical systems, see \cite{Balakrishnan1960},\cite{Komatsu1966},\cite{HA}, \cite{CMCMA}. The framework of abstract parabolic problems with critical nonlinearities \cite{JAAC} further underscores the importance of fractional calculus in PDE analysis. 
  Recent extensions to non-commutative settings, particularly for quaternionic linear operators \cite{ColomboGantner2019}, have revealed new challenges in spectral analysis and functional calculus due to the non-commutative nature of quaternionic multiplication. By extending this framework to block operator matrices, 
  we aim to provide a deeper understanding of the interplay between operator theory and the analysis of coupled PDE systems.
Block operator matrices are widely used in the modeling of multi-component systems where multiple 
state variables and their interactions must be considered simultaneously. These matrices appear naturally in applications such as reaction-diffusion systems, evolution equations, quantum mechanics, control theory, and magnetohydrodynamics, see \cite{KJE}, \cite{FAHARMS}, \cite{RNA}, with recent advances in Schr{\"o}dinger cascade systems \cite{MNK} demonstrating their fractional power behavior. 
In quaternionic settings, such matrices model systems with intrinsic non-commutative symmetries, though their fractional powers remain largely unexplored. The spectral challenges parallel those in unbounded operator matrices \cite{RN}, where non-commutativity complicates resolvent analysis. In such settings, 
fractional powers of operators serve not only as a generalization of integer powers but also as essential instruments
 to characterize the behavior of systems governed by non-local phenomena, anomalous diffusion, and fractional time derivatives \cite{BCCN}, \cite{BCN}, extending the infinite-dimensional dynamical systems framework \cite{JCR}.
Inspired by \cite{Belluzi}, our primary objective in this work, is to build a rigorous and comprehensive theoretical framework for defining and calculating 
the fractional powers of positive operators, and subsequently extend these results to  
$3\times3$ block operator matrices. A fundamental open question emerges: Can our methods be extended to quaternionic block operator matrices, where the non-commutativity of entries introduces additional complexity in defining consistent fractional powers? Positive operators, which are central to our approach, 
exhibit spectral properties that facilitate the definition of fractional powers through complex function theory. 
The extension to block operator matrices introduces additional challenges due to the potential non-commutativity of the operator entries, 
the unboundedness of the components, and the intricate structure of their spectra, see \cite{AJER}, \cite{SCAJRM}, \cite{RN2}, \cite{CT}. Semigroup generation techniques \cite{JLJHAC} offer complementary approaches to these challenges. 
To address these challenges, we employ several advanced analytical techniques to derive explicit expressions for the fractional powers of block operator matrices. First, we utilize alternative formulas grounded in spectral theory, 
 which allow for precise characterization of fractional powers based on the resolvent operators. Next, we apply changes of variables to transform complex integrals into more tractable forms, 
 facilitating the evaluation of fractional powers in specific cases. Furthermore, we employ the second resolvent identity, 
 a powerful tool in perturbation theory, to handle cases where the block components differ significantly.
Each of these methods is carefully analyzed to identify their advantages,
 limitations, and applicability to various types of block operator matrices. 
 Our theoretical findings are then applied to concrete examples of PDE systems, 
 demonstrating the effectiveness of fractional operator theory in capturing the nuanced dynamics of multi-component systems. These applications include reaction-diffusion systems, coupled oscillators, 
  and diffusion-advection equations, where the use of fractional powers enhances the modeling capabilities by accounting for memory effects, long-range interactions, and non-local phenomena, mirroring geometric theories for semilinear equations \cite{DH}.
The structure of the paper is organized as follows:
 we begin with a detailed review of the mathematical preliminaries necessary for understanding positive operators, 
 their spectra, and fractional powers. 
 This includes a discussion on the holomorphic functional calculus, semi-group theory, 
 and the resolvent approach. Subsequently, we present the main theoretical results on the fractional powers of 
$3\times3$ block operator matrices, emphasizing different computational techniques and the conditions under which they apply.
 We then provide a series of illustrative examples drawn from PDEs to showcase the practical relevance of our methods. Finally, 
 we discuss the implications of our findings for future research in operator theory, particularly in the context of higher-dimensional 
 block matrices and more complex dynamical systems.
By bridging the gap between abstract operator theory and applied PDE analysis,
 this work offers significant insights into the use of fractional powers for block operator matrices. 
 It not only enriches the theoretical understanding of operator powers but also opens up new avenues
  for the modeling and analysis of complex multi-component systems.
   The results presented here have the potential to inspire further research in fractional calculus, 
   functional analysis, and the study of advanced PDE systems.
\section{Mathematical preliminaries}
In this section, we present general definitions, results and notaions about positive operators. 
In fact, we investigate the definition of fractional powers for these positive operators and introduce alternative approaches that simplify the process of calculating fractional powers.
\begin{Definition}[\cite{HA}]
Let $X$ be a Banach space and $ A:D(A)\subset X \rightarrow X$ a closed operator with dense domain. We say that $A$ is positive, if there exists $M\geq1$ such that:
\begin{enumerate}
\item $[0,\infty)\subset\rho(-A),$
\item $(1+s)\|(s+A)^{-1}\|_{\mathcal{L} (X)}\leq M$, for any positive real $s$.
\end{enumerate}
 We denote $\mathcal{P}$, the class of positive operators on $X$. If we want to emphasize the constant $M$, we denote $A\in\mathcal{P}_M(X)$.
\end{Definition}

\begin{Proposition}[\cite{HA}, Section 4.9]
Let $A\in \mathcal{P}_M(X)$. If $\theta_M=\arcsin(\frac{1}{2M})$, then the set
\begin{align*} \Xi_{M}=\Big\{z\in\mathbb{C}:|\arg z|\leq\theta_{M}\Big\}\cup\Big\{z\in\mathbb{C}:|z|\leq\frac{1}{2M}\Big\} \end{align*}
is is contained in $\rho(-A)$ and we have the following estimate
\begin{equation} (1+|\lambda|)\|(\lambda+A)^{-1}\|_{\mathcal{L} (X)}\leq2M+1,\ \forall\lambda\in\Xi_{M}.
\label{estim}
\end{equation}
\end{Proposition}

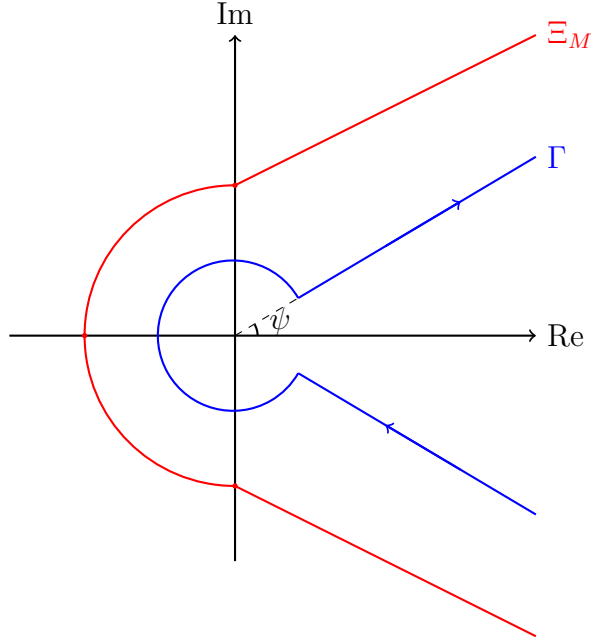
\begin{figure}[H]
    \centering
    \begin{tikzpicture}[scale=2, domain=0:2]

        \draw[->,thick] (-1.5,0) -- (2,0) node[right] {$\mathrm{Re}$}; 
        \draw[->,thick] (0,-1.5) -- (0,2) node[above] {$\mathrm{Im}$}; 

        \draw[thick,red,domain=0:2] plot (\x, {0.5*\x + 1}) node[right] {$\Xi_{M}$};

        \draw[thick,red,domain=0:2] plot (\x, {-0.5*\x - 1});

        \draw[thick,red] (0,1) arc (90:270:1cm);

        \draw[thick,blue] (0.42,0.25) arc (30:330:0.5cm);

        \draw[thick,blue,domain=0.42:2] plot (\x, {0.595*\x}) node[right] {$\Gamma$};

        \draw[thick,blue,domain=0.42:2] plot (\x, {-0.595*\x});

        \draw[thick,blue,->] (1, {0.595*1}) -- (1.5, {0.595*1.5});

        \draw[thick,blue,<-] (1, {-0.595*1}) -- (1.5, {-0.595*1.5});

        \fill[red] (0,1) circle (0.5pt); 
        \fill[red] (-1,0) circle (0.5pt); 
        \fill[red] (0,-1) circle (0.5pt); 

        \draw[dashed] (0,0) -- (0.42,0.25);

        \draw[thick] (0.15, 0) arc (0:30:0.15);
        \node at (0.3, 0.1) {$\psi$}; 

    \end{tikzpicture}
    \caption{Sets in $\rho(-A)$}
    \label{figureone}
\end{figure}

 Since $\Xi_{M}\subset\rho(-A)$ (see Figure \ref{figureone}), we can find $\psi\in(0,\theta_{M})$ and $r_{0}\in(0,\frac{1}{2M})$
such that if $\Sigma_{\psi}=\{z\in\mathbb{C}:|\arg z|\leq\psi\}$ then \begin{align*}\Sigma_{\psi}\cup\mathrm{B}[0,r_{0}]\subset\rho(-A).\end{align*}
 This figure can be found in \cite{Belluzi}.

The fractional powers of a positive operator $A$ (in the sense of \cite{HA}) are well defined and given in the following. Let us denote by $\Phi$ the set of holomorphic functions $ \varphi\in\mathcal{H}(\mathbb{C\setminus R-})$ such that: there exists $\delta_{\varphi}>0$ with

\begin{equation}\label{amann466} |\lambda|^{\delta_{\varphi}}\varphi(\lambda)\rightarrow 0 \,\,\mbox{when}\,\, |\lambda|\rightarrow\infty
\,\mbox{uniformly in}\, \{|\arg(\lambda)|\leq \pi-\epsilon\}\, \mbox{for all}\,\epsilon\in (0,\pi).\end{equation}
$\Phi$ is a commutative algebra without unit (for point-by-point multiplication).
Let $\varphi\in\Phi$, note:
\begin{equation}\label{ammann467} \varphi(A):=\frac{1}{2\pi i}\int_{\Gamma}\varphi(-\lambda)(\lambda+A)^{-1}d\lambda=\frac{1}{2\pi i}\int_{-\Gamma}\varphi(\lambda)(\lambda-A)^{-1}d\lambda,\end{equation}
where $\Gamma$ is a contour in $\Sigma_{\psi}\cup\mathrm{B}(0,r_{0})-\mathbb{R_{+}}$ oriented from $\infty e^{-i\psi}$ to $\infty e^{i\psi}$ for certain $\psi\in(0,\arcsin(\frac{1}{2M}))$.
From (\ref{amann466}), (\ref{estim}) and the Cauchy's Integral Theorem, it follows that $\varphi(A)$ is well-defined, independent of the choice of $\Gamma$ and $\varphi(A)\in \mathcal{L}(X)$. Furthermore, we have :
\begin{Lemma}\cite[Lemma 4.6.1]{HA}
\begin{equation*} \Phi\rightarrow \mathcal{L}(X)\,,\varphi\mapsto\varphi(A)\end{equation*}
is a homomorphism of algebras.\label{morphalg}
\end{Lemma}
Let $z\in\mathbb{C}$ and $\varphi_{z}(\lambda):=\lambda^{z}=e^{z\log(\lambda)},\lambda\in\mathbb{C\setminus R-}$, where $\log$ is the principal determination of the logarithm. Note that $\varphi_{z}\in\Phi$ for $Re(z)<0$.
If $k\in\mathbb{N^{*}}$ then $\varphi_{-k}\in\mathcal{H}(\mathbb{C^{*}})$ and we can deform $\Gamma$ to obtain a positively oriented circle of center $0$ and radius $r$, with $0<r<\frac{1}{\|A^{-1}\|}.$ Since \begin{align*}(\lambda-A)^{-1}=-A^{-1}(1-\lambda A^{-1})^{-1}=-\sum_{j=0}^{\infty}A^{-j-1}\lambda^{j},\,\, |\lambda|=r,\end{align*} it follows that:
\begin{align} \varphi_{-k}(A)=\frac{-1}{2\pi i}\int_{|\lambda|=r}\lambda^{-k}(\lambda-A)^{-1}d\lambda=\sum_{j=0}^{\infty}A^{-j-1}Res(\lambda^{j-k},0)=A^{-k}.\label{4.6.9}\end{align}
This justifies the following definition.
\vskip 1cm
The fractional power of a positive operator $A$ of exponent $z\in\mathbb{C}$ with $Re(z)<0$
denoted $A^{z}$ is equal to $\varphi_{z}(A)$.
\begin{Definition}[\cite{HA}]
Let $A\in \mathcal{P}_M(X)$ and $\Gamma$ be a $\Sigma_{\psi}\cup\mathrm{B}(0,r_{0})-\mathbb{R_{+}}$ contour oriented from $\infty e^{-i\psi}$ to $\infty e^{i\psi}$. For all $\alpha>0$, we define:
 \begin{align} A^{-\alpha}:=\frac{1}{2i\pi}\int_{\Gamma}(-\lambda)^{-\alpha}(\lambda+A)^{-1}d\lambda.
 \end{align}
\end{Definition}
\begin{Remark}
In the context of holomorphic functional calculus for a bounded operator $T$, we know that $\sigma(T)$ is compact, so we can surround it with a simple, closed contour. When manipulating this type of integral with a closed, unbounded operator, we lose the boundedness of the spectrum, so the contour chosen is expected to be infinite.
\end{Remark}
\begin{Definition} A family $\{T(t), t\geq0\}\subset\mathcal{L}(X)$ is said to be \emph{semi-group of bounded operators} on $X$ if:
\begin{itemize}
  \item $T(0)=I_{X}$
  \item $T(t+s)=T(t)T(s),\forall t,s\in[0,+\infty)$
\end{itemize}
\end{Definition}
\begin{Definition}\cite{} A semi-group $\{T(t), t\geq0\}\subset\mathcal{L}(X)$ is a \textup{strongly continuous} semi-group ($\mathcal{C}_{0}$-semi-group), if:
\begin{align*}\lim\limits_{t\rightarrow 0} T(t)x=x,\forall x\in X.\end{align*}
The linear operator: \begin{align*}A : D(A)\subset X\rightarrow X \end{align*}
defined by \begin{align*} D(A)=\Big\{x\in X;\,\lim\limits_{t\rightarrow 0^{+}}\frac{T(t)x-x}{t}\,exists\Big\} \end{align*}
and
\begin{align*}
Ax = \lim\limits_{t\rightarrow 0^{+}} \frac{T(t)x - x}{t}, \quad \forall x\in D(A)
\end{align*}
is called \textup{infinitesimal generator} of the $\mathcal{C}_{0}$-semi-group of $\{T(t), t\geq0\}$.
\end{Definition}
\begin{Lemma}\label{unifbor} Let $A\in\mathcal{P}_{M}(X)$, the family $(A^{-\alpha})_{\alpha\in(0,1)}$ is uniformly bounded\end{Lemma}

\proof
Let's fix $\alpha\in(0,1)$.
\begin{align*}A^{-\alpha}&=\frac{1}{2i\pi}\int_{\Gamma}(-\lambda)^{-\alpha}(\lambda+A)^{-1}d\lambda\\
&=\frac{e^{i\pi\alpha}}{2i\pi}\int_{0}^{\infty}s^{-\alpha}(s+A)^{-1}ds\\
&-\frac{e^{-i\pi\alpha}}{2i\pi}\int_{0}^{\infty}s^{-\alpha}(s+A)^{-1}ds\\
&=\frac{\sin(\pi\alpha)}{\pi}\int_{0}^{\infty}s^{-\alpha}(s+A)^{-1}ds.
\end{align*}
Thus, \begin{align}\label{4.6.9} A^{-\alpha}=\frac{\sin(\pi\alpha)}{\pi}\int_{0}^{\infty}s^{-\alpha}(s+A)^{-1}ds,\,\forall\alpha\in(0,1)\end{align}
In particular, if $X=\mathbb{C},\,\mbox{and}\,A=1$ we find \begin{align}\int_{0}^{\infty}s^{-\alpha}(s+A)^{-1}ds=\frac{\pi}{\sin(\pi\alpha)},\,\,0<\alpha<1.\label{4.6.10}\end{align}
We deduce from (\ref{4.6.9}) and (\ref{4.6.10}) that \begin{align*} \|A^{-\alpha}||\leq M,\,\forall\alpha\in (0,1).\end{align*}
\qed
\begin{theorem}\cite[Theorem 4.6.2]{HA} For any positive operator $A\in\mathcal{P}_{M}(X)$, the family $(A^{-\alpha})_{\alpha\geq0}$ is a $\mathcal{C}_{0}$-semi-group.
Let $-\log A$ be its infinitesimal generator.
\end{theorem}

\begin{Proposition}\label{injective}\cite[Proposition 5.30]{EKNR} For any $\alpha>0$, the operator $A^{-\alpha}\in\mathcal{L}(X)$ is injective.
\end{Proposition}

\begin{Definition} For $\alpha>0$ we define $A^{\alpha}$, as the inverse of $A^{-\alpha}$. In other words:
\begin{align*} A^{\alpha}:=(A^{-\alpha})^{-1}. \end{align*}
 In this case, $D(A^{\alpha})=R(A^{-\alpha})$ and we denote $X^{\alpha}$ the domain of $A^{\alpha}$ equipped with the graph norm: $X^{\alpha}=(D(A^{\alpha}),\|.\|_{X^{\alpha}})$ and:
\begin{align*} \|u\|_{X^{\alpha}}=\|u\|+\|A^{\alpha}u\|, \forall u\in D(A^{\alpha}).\end{align*}.
\end{Definition}
\begin{Proposition}
For a positive operator $A$, we have the following properties:
\begin{enumerate}
  \item\label{p1}For $\alpha>0$, $A^{\alpha}$ is a closed operator.
  \item\label{p2} Let $\alpha,\beta>0$ with $\alpha<\beta$, then $D(A^{\beta})\subset D(A^{\alpha})$.
  \item\label{p3} $\overline{D(A^{\alpha})}=X$, $\forall\alpha>0$.

\end{enumerate}
\label{normban}\end{Proposition}
\proof

(\ref{p1}): Let $(x_{n})\in D(A^{\alpha})^{\mathbb{N}}$ be such that $x_{n}\rightarrow x$ and $A^{\alpha}x_{n}\rightarrow y$.
$A^{-\alpha}$ is bounded ( and therefore closed), from $x_{n}\rightarrow A^{-\alpha}y$ and $A^{-\alpha}x_{n}\rightarrow A^{-\alpha}x$.
 We find $x=A^{-\alpha}y\in R(A^{-\alpha})=D(A^{-\alpha})$ and $A^{-\alpha}x=y$. Therefore, $A^{\alpha}$ is a closed operator.\\ (\ref{p2}): Let $x\in D(A^{\beta})$, we have $x=A^{-\alpha}(A^{-\beta+\alpha}A^{\beta}x)$ so $x\in R(A^{-\alpha})=D(A^{\alpha})$.\\
(\ref{p3}): Let $x\in D(A)$ and $y=Ax$. We know that $\overline{D(A)}=X$. Let's set $\epsilon>0$.
\begin{align*} \exists u\in\ D(A);\,\|u-y\|\leq\frac{\epsilon}{\|A^{-1}\|}.\end{align*}
Let's put $v=Au$. We have:
\begin{align*} \|A^{-2}v-x\|=\|A^{-1}u-A^{-1}y\|\leq\|A^{-1}\| \|u-y\|\leq\epsilon.\end{align*}
Thus, $\overline{D(A^{2})}\supset D(A)$ or $\overline{D(A^2)}\supset \overline{D(A)}=X$.
By induction on $k\in\mathbb{N}$ we have $\overline{D(A^{k})}=X$.
Let $\alpha>0$ and $k\in\mathbb{N}$ be such that $\alpha<k$.
From (\ref{p2}), we obtain
\begin{align*}D(A^{k})\subset D(A^{\alpha}).\end{align*}
 Or else:
\begin{align*}X=\overline{D(A^{k})}\subset \overline{D(A^{\alpha})}.\end{align*}
\qed
\begin{Remark} Propiety (1) of (Proposition \ref{normban}) is equivalent to: $(D(A^{\alpha}),\|\,\|_{X^{\alpha}})$ is a Banach space.
\end{Remark}
\begin{Proposition}
Let $\alpha>0$ and $A\in\mathcal{P}(X)$.
$\|.\|_{\alpha}\, : \, x\ni D(A^{\alpha})\rightarrow \|A^{\alpha}x\|$ is a norm on $D(A^{\alpha})$.
Moreover, $(D(A^{\alpha}), \|.\|_{\alpha})$ is a Banach space and $\|.\|_{\alpha}$ is equivalent to $\|.\|_{X^{\alpha}}$.
\end{Proposition}
\proof From the injectivity of $A^{\alpha}$ and $\|.\|$ is a norm on $X$, we deduce that $\|.\|_{\alpha}$ is a norm on $D(A^{\alpha})$. Let $x_{n}$ be a Cauchy sequence in $(D(A^{\alpha}), \|.\|_{\alpha})$, then $A^{\alpha}x_{n}$ is Cauchy in the complete $(X, \|.\|)$ and converges to a certain
$y\in X$. The continuity of $A^{-\alpha}$ implies that $x_{n}$ converges in $X$ to $A^{-\alpha}y$. We have $A^{-\alpha}y\in D(A^{\alpha})$ and $A^{\alpha}A^{-\alpha}y=y$,
 which gives $\|x_{n}-A^{-\alpha}y\|_{\alpha}\rightarrow 0$, when $n\rightarrow +\infty$. $(D(A^{\alpha}), \|.\|_{\alpha})$ is therefore Banach space.
 The two norms $ \|.\|_{\alpha}$ and $ \|.\|_{X^{\alpha}}$ are comparable, Banach's isomorphism theorem ensures that they are equivalent.\qed
\begin{Proposition}
If $0<\alpha<\beta,$ then the injection $X^{\beta}\hookrightarrow X^{\alpha}$ is continuous.
\end{Proposition}
\proof From (Proposition \ref{normban}), we have $X^{\beta}\subset X^{\alpha}$. If $x\in X^{\beta}$, then:
\begin{align*}\|x\|_{\alpha}=\|A^{\alpha}x\|=\|A^{\alpha-\beta}A^{\beta}x\|& \leq\|A^{\alpha-\beta}\| \|A^{\beta}x\|\\
&\leq \|A^{\alpha-\beta}\| \|x\|_{\beta}.
\end{align*}
\qed
\begin{Proposition} \label{prop25} \cite[(4.6.9), Theorems 4.6.3 and 4.6.5]{HA}
Let $A\in \mathcal{P}(X)$:
\begin{enumerate}
    \item If $\alpha\in(0,1)$, then:
    \begin{align}\label{e.1}
        A^{-\alpha} = \frac{\sin(\pi\alpha)}{\pi} \int_{0}^{\infty} s^{-\alpha} (s+A)^{-1} ds.
    \end{align}

    \item If $\alpha\in(0,m+1)$, $\alpha$ non-integer, $m\in\mathbb{N}^*$, then:
    \begin{align}\label{e.2}
        A^{-\alpha} = \frac{\sin(\pi\alpha)}{\pi} \frac{m!}{(1-\alpha)(2-\alpha)\dots(m-\alpha)}
        \int_{0}^{\infty} s^{m-\alpha} (s+A)^{-m-1} ds.
    \end{align}

    \item If $-1<\alpha<1$ and $x\in D(A)$, then:
    \begin{align}\label{e.3}
        A^{-\alpha}x = \frac{\sin(\pi\alpha)}{\pi\alpha} \int_{0}^{\infty} s^{\alpha} (s+A)^{-2} Ax \, ds.
    \end{align}
\end{enumerate}
\end{Proposition}
If $A$ is positive, the spectrum of $A^\alpha$, denoted by $\sigma(A^\alpha)$, can be fully determined from the spectrum of $A$, $\sigma(A)$. The relationship between these spectra is straightforward, as described in the following proposition.
\begin{Proposition}\cite{CMCMA} If $A \in \mathcal{P}(X)$, then for each $\alpha > 0$, we have:
\begin{align*}
\sigma(A^{\alpha}) = [\sigma(A)]^{\alpha} = \left\{\sigma^{\alpha} \mid \sigma \in \sigma(A)\right\}
\end{align*}
\end{Proposition}
\section{Fractional powers of $3\times3$ block operators matrix}
 In this section, we'll give examples of $3\times3$ matrix operators whose fractional powers can be calculated.
In the following, $X$ is a Banach space and $A$ is a positive operator in $X$.
\\In many applications, it is typical to handle systems of partial differential equations rather than dealing with a single equation.
 In these cases, the equation is generally put into its matrix formulation, and the linear operator can be given by a matrix operator whose entrances are unbounded linear operators.
In an abstract general framework, let's consider $\Lambda : D(\Lambda)\subset X_{1}\times X_{2}\times X_{3}\rightarrow X_{1}\times X_{2}\times X_{3}$
 a matrix whose entries are linear operators that commute in an appropriate domain, where $X_{1}$, $X_{2}$ $X_{3}$ are Banach spaces (which are in general sub-vector spaces of $X$ with certain inclusion and density properties).
 We note:
 \begin{align*}
 \Lambda=\left(
 \begin{array}{ccc}
  A_{11} & A_{12} & A_{13} \\
  A_{21} & A_{22} & A_{23} \\
  A_{31} & A_{32} & A_{33} \\
 \end{array}
 \right)\end{align*}
with $A_{ij} : D(\Lambda_{ij})\subset X_{j}\rightarrow X_{i}$, $i,j\in \{1,2,3\}$.\\
In the following, we assume that:\\
\begin{enumerate}
\item $\Lambda$ is a positive operator. From now on, for any $0<\alpha<1$, its fractional power is given by:
\begin{align*} \Lambda^{-\alpha}=\frac{\sin(\pi\alpha)}{\pi}\int_{0}^{\infty}s^{-\alpha}(s+\Lambda)^{-1}ds\end{align*}\label{cfh}
\item For all $s>0$, $(s+\Lambda)^{-1}$ is given by:
\begin{equation*}
\left(
\begin{array}{ccc}
(s+A_{22})(s+A_{33})-A_{23}A_{32} & A_{32}A_{13}-A_{12}(s+A_{33}) & A_{12}A_{23}-A_{13}(s+A_{22})\\
A_{13}A_{23}-A_{21}(s+A_{33}) & (s+A_{11})(s+A_{33})-A_{13}A_{31} & A_{21}A_{13}-A_{23}(s+A_{11})\\
A_{21}A_{32}-A_{31}(s+A_{22}) & A_{31}A_{12}-A_{32}(s+A_{11}) & (s+A_{22})(s+A_{11})-A_{21}A_{12}\\
\end{array}\right)[det(\Lambda)]^{-1}.
\label{hyp}
\end{equation*}
\end{enumerate}
 with
\begin{align*}
[det(\Lambda)]^{-1}
&=[(s+A_{11})(s+A_{22})(s+A_{33})-A_{21}A_{12}(s+A_{33})-A_{31}A_{13}(s+A_{22})\\
&-A_{23}A_{32}(s+A_{11})+A_{21}A_{13}A_{32}+A_{31}A_{12}A_{23}]^{-1}.
\end{align*}
 Simply put, we can treat the inverse of $(s+\Lambda)$ as if it were a $3 \times 3$ matrix with complex entries. Nevertheless, it's important to note that equality does not generally hold, so our focus will be on cases where equality does hold.
To simplify \eqref{hyp}, we observe:
   \begin{align*}
     (s+\Lambda)^{-1}=\left(
                                                                        \begin{array}{ccc}
                                                                          \widetilde{A}_{11}(s) & \widetilde{A}_{12}(s) & \widetilde{A}_{13}(s) \\
                                                                          \widetilde{A}_{21}(s) & \widetilde{A}_{22}(s) & \widetilde{A}_{23}(s) \\
                                                                          \widetilde{A}_{31}(s) & \widetilde{A}_{32}(s) & \widetilde{A}_{33}(s) \\
                                                                        \end{array}
                                                                      \right).\label{tild}\end{align*}
Where $\widetilde{A}_{ij}(s)=(s+\Lambda)_{ij}^{-1}[det(\Lambda)]^{-1}, \;\forall i,j\in\{1,2,3\}$.
In this case, expression\eqref{cfh} gives:
\begin{equation}\Lambda^{-\alpha}=\left(
 \begin{array}{ccc}
  \frac{\sin(\pi\alpha)}{\pi}\int_{0}^{\infty}s^{-\alpha}\widetilde{A}_{11}(s)ds & \frac{\sin(\pi\alpha)}{\pi}\int_{0}^{\infty}s^{-\alpha}\widetilde{A}_{12}(s)ds & \frac{\sin(\pi\alpha)}{\pi}\int_{0}^{\infty}s^{-\alpha}\widetilde{A}_{13}(s)ds \\
  \frac{\sin(\pi\alpha)}{\pi}\int_{0}^{\infty}s^{-\alpha}\widetilde{A}_{21}(s)ds & \frac{\sin(\pi\alpha)}{\pi}\int_{0}^{\infty}s^{-\alpha}\widetilde{A}_{22}(s)ds& \frac{\sin(\pi\alpha)}{\pi}\int_{0}^{\infty}s^{-\alpha}\widetilde{A}_{23}(s)ds \\
  \frac{\sin(\pi\alpha)}{\pi}\int_{0}^{\infty}s^{-\alpha}\widetilde{A}_{31}(s)ds & \frac{\sin(\pi\alpha)}{\pi}\int_{0}^{\infty}s^{-\alpha}\widetilde{A}_{32}(s)ds & \frac{\sin(\pi\alpha)}{\pi}\int_{0}^{\infty}s^{-\alpha}\widetilde{A}_{33}(s)ds \\
  \end{array}
  \right).
  \label{puiss-alpha}
\end{equation}
\subsection{Using alternative formulas}
The simplest situation would occur if $\widetilde{\Lambda}_{ij}(s)$ took a form that allowed us to apply the characterization from Proposition \eqref{prop25}. For example, if $\widetilde{\Lambda}_{ij}(s)= (s + A)^{-1}$, where $A$ is a positive operator, we would have:

\begin{align*}
\frac{\sin(\pi\alpha)}{\pi}
\int_0^\infty s^{-\alpha} \widetilde{\Lambda}_{ij}(s) \, ds = \frac{\sin(\pi\alpha)}{\pi}
\int_0^\infty s^{-\alpha} (s + A)^{-1} \, ds = A^{-\alpha}, \quad \text{for } 0 < \alpha < 1.
\end{align*}

\begin{example}We consider the following DPE, where $\Omega$ is a smooth boundary domain and $-\Delta_{D}$ is the positive Laplacian with Dirichlet boundary conditions.
\begin{equation}
    \begin{cases}
      u_{t}-\Delta_{D}u=f, \;  x\in\Omega, t>0,\\
      v_{t}-\Delta_{D}v=g,  \;x\in\Omega, t>0,\\
      w_{t}+u-\Delta_{D}w=g,\;  x\in\Omega, t>0,\\
      u(x,0)=v(x,0)=w(x,0)=0,\;             x\in\partial\Omega.\\
    \end{cases}\,.
\label{edp1}
\end{equation}
                     \end{example}
 Let's consider the matrix operator:\\
  $\Lambda_{1}: D(A)\times D(A)\times D(A)\subset X\times X\times X\rightarrow X\times X\times X$ given by:
  \\
 \begin{equation}\Lambda_{1}= \left(
                   \begin{array}{ccc}
                     A & 0 & 0 \\
                     0 & A & 0 \\
                     I & 0 & A \\
                   \end{array}
                 \right)
 \label{first}
 \end{equation} where $A=-\Delta_{D}$.
 \begin{Lemma} If $ A : D(A)\subset X\rightarrow X $ is positive and $0<\gamma<2$,\; $\gamma\neq1$,\; then:
 \begin{align*} \int_{0}^{\infty} s^{1-\gamma}(s+A)^{-2}ds=\frac{\pi}{\sin(\pi\gamma)}(1-\gamma)A^{-\gamma}.
 \end{align*}
 \label{lemme1}
\end{Lemma}
\proof
\begin{align*}
\int_{0}^{\infty}s^{1-\gamma}(s+A)^{-2} \, ds 
&= \frac{\pi}{\sin(\pi\gamma)}(1-\gamma) \frac{\sin(\pi\gamma)}{\pi} \frac{1}{1-\gamma} \int_{0}^{\infty} s^{1-\gamma}(s+A)^{-2} \, ds \\
&= \frac{\pi}{\sin(\pi\gamma)}(1-\gamma) A^{-\gamma}.\end{align*}
\qed\\
In particular, for $-1<\alpha<1$, $\alpha\neq0$ we have:
\begin{align*}
\frac{\sin(\pi\alpha)}{\pi}\int_{0}^{\infty}s^{-\alpha}(s+A)^{-2}ds&=\frac{\sin(\pi\alpha)}{\pi}\int_{0}^{\infty}s^{1-(1+\alpha)}(s+A)^{-2}ds\\&= \frac{\sin(\pi\alpha)}{\pi}\frac{\pi}{\sin(\alpha\pi+\pi)}(1-(1+\alpha))A^{-1-\alpha}\\&=\alpha A^{-1-\alpha}.\end{align*}
 \begin{Lemma}$\Lambda_{1}$ with $D(\Lambda_{1})=X^{1}\times X^{1}\times X^{1}$ is closed with dense domaine.\end{Lemma}
\proof Since $\overline{X^{1}}=X$, then $\overline{D\left(\Lambda_{1}\right)}=X \times X \times X$. Let

\begin{align*}
U_{n}=\left(\begin{array}{c}
u_{n} \\
v_{n} \\
w_{n}
\end{array}\right) \in D\left(\Lambda_{1}\right)^{\mathbb{N}}
\end{align*}

such that $U_{n} \rightarrow U=\left(\begin{array}{c}u \\ v \\ w\end{array}\right)$ and $\Lambda_{1} U_{n} \rightarrow V=\left(\begin{array}{l}x \\ y \\ z\end{array}\right)$. We have

\begin{align*}
\left\{\begin{array}{l}
u_{n} \rightarrow u \quad \text { and } \quad A u_{n} \rightarrow x \\
v_{n} \rightarrow v \quad \text { and } A v_{n} \rightarrow y \\
w_{n} \rightarrow w \quad \text { and } A w_{n} \rightarrow z-x
\end{array}\right.
\end{align*}
Since $A$ is closed, then $u, v, w \in D(A)=X^{1}$ and

\begin{align*}
\left\{\begin{array}{l}
A u=x \\
A v=y \\
x-A w=z
\end{array}\right.
\end{align*}

Thus, $U \in D\left(\Lambda_{1}\right)$ and $\Lambda_{1} U=V$. Therefore, $\Lambda_{1}$ is closed.
   \qed
 \begin{Proposition}
    If $ A : D(A)\subset X\rightarrow X $ is a positive operator, and $\Lambda_{1}$ is the operator given in \eqref{first}, then:
 \begin{enumerate}
 \item $\Lambda_{1}$ is a positive operator and $D(\Lambda_{1})=X^{1}\times X^{1}\times X^{1}$\label{p1}
 \item $\sigma(\Lambda_{1})=\sigma(A)$\label{p2}
 \item $\forall \lambda\in\rho(-\Lambda_{1})$\label{p3}, we have:
 \begin{align*}
 (\lambda+\Lambda)^{-1}=\left(
                          \begin{array}{ccc}
                            (\lambda+A)^{-1} & 0 & 0 \\
                            0 & (\lambda+A)^{-1} & 0 \\
                            -(\lambda+A)^{-2} & 0 & (\lambda+A)^{-1} \\
                          \end{array}
                        \right).
 \end{align*}
 \end{enumerate}
 \label{prop1}
 \end{Proposition}
 \proof \eqref{p1} Let $s\in[0, \infty)$, and $U=\left(\begin{array}{c}u \\ v \\ w\end{array}\right) \in D\left(\Lambda_{1}\right)$ such that

\begin{align*}
\left(s+\Lambda_{1}\right) U=\left(\begin{array}{l}
0 \\
0 \\
0
\end{array}\right).
\end{align*}
We have

\begin{align*}
\left\{\begin{array}{l}
(s+A) u=0 \\
(s+A) v=0 \\
u+(s+A) w=0
\end{array}\right.
\end{align*}

Since $A$ is positive operator then $u=v=w=0$, thus $s+\Lambda_{1}$ is injective. Set \begin{align*} V=\left(\begin{array}{l}x \\ y \\ z\end{array}\right) \in X \times X \times X \,\mbox{and}\, U=\left(\begin{array}{c}(s+A)^{-1} x \\ (s+A)^{-1} y \\ x+(s+A)^{-1}(z-x)\end{array}\right).\end{align*} We have $\left(s+\lambda_{1}\right) U=V$, so $s+\Lambda_{1}: D\left(\Lambda_{1}\right) \subset X \times X \times X \rightarrow X \times X \times X$ is surjective and therefore bijective. Since $s+\Lambda_{1}$ is closed, then, by the closed graph theorem, $\left(s+\Lambda_{1}\right)^{-1}$ is bounded. Thus, $[0, \infty) \subset \rho\left(-\Lambda_{1}\right)$, we can verify that

\begin{align*}
\left(s+\Lambda_{1}\right)^{-1}=\left(\begin{array}{ccc}
(s+A)^{-1} & 0 & 0 \\
0 & (s+A)^{-1} & 0 \\
-(s+A)^{-2} & 0 & (s+A)^{-1}
\end{array}\right)
\end{align*}
and
\begin{align*}
(1+s)\left\|\left(s+\Lambda_{1}\right)^{-1}\right\|_{\mathcal{L}(X \times X \times X)} &\leq 3(1+s)\left\|(s+A)^{-1}\right\|_{\mathcal{L}(X)} \\
&\quad + (1+s)\left\|(s+A)^{-1}\right\|_{\mathcal{L}(X)}\left\|(s+A)^{-1}\right\|_{\mathcal{L}(X)} \\
&\leq 3 M + \frac{M^{2}}{1+s} \\
&\leq 3 M + M^{2}.
\end{align*}
This implies that $\Lambda_{1}$ is a positive operator.\\
\eqref{p2} Let $\lambda \in \rho(-A)$, by adopting the same approach used to show that $[0, \infty) \subset \rho\left(-\Lambda_{1}\right)$ we have $\lambda+\Lambda_{1}: D\left(\Lambda_{1}\right) \rightarrow X \times X \times X$ is bijective, thus $\rho(A) \subset \rho\left(\Lambda_{1}\right)$. Take $\lambda \in \rho\left(-\Lambda_{1}\right)$, if for some $u \in D(A)$ we have $(\lambda+A) u=0$, then \begin{align*}\left(\lambda+\Lambda_{1}\right)\left(\begin{array}{l}u \\ 0 \\ 0\end{array}\right)=\left(\begin{array}{l}0 \\ 0 \\ 0\end{array}\right),\end{align*} thus $u=0, \lambda+A$ is injective. For the surjectivity, let $y \in X$ and $Y=\left(\begin{array}{l}y \\ 0 \\ 0\end{array}\right)$. Since $\lambda+\Lambda_{1}$ is surjective, there exists $U=\left(\begin{array}{c}u \\ v \\ w\end{array}\right) \in D(A) \times D(A) \times D(A)$ such that $\left(\lambda+\Lambda_{1}\right) U=Y$. Hence $(\lambda+A) u=y,$ by closed graph theorem, one has $\lambda \in \rho(-A)$. Thus $\rho\left(\Lambda_{1}\right)=\rho(A)$.\\ \eqref{p3} easy to check.
\qed
\begin{theorem}
  Let $ A : D(A)\subset X\rightarrow X $ be a positive operator. \\ The Fractional powers of $\Lambda_{1}$ given by \eqref{first}, can be calculated explicitly.\\
  For all $\alpha\in(0,1)$, we have:
  \begin{align*}
 \Lambda_{1}^{-\alpha}=\left(
                          \begin{array}{ccc}
                            A^{-\alpha} & 0 & 0 \\
                            0 & A^{-\alpha} & 0 \\
                            -\alpha A^{-\alpha-1} & 0 & A^{-\alpha} \\
                          \end{array}
                        \right),\; with\; \Lambda_{1}^{-\alpha}\in\mathcal{L}(X\times X\times X),
 \end{align*}\\and:
 \begin{align*}
 \Lambda_{1}^{\alpha}=\left(
                          \begin{array}{ccc}
                            A^{\alpha} & 0 & 0 \\
                            0 & A^{\alpha} & 0 \\
                            \alpha A^{\alpha-1} & 0 & A^{\alpha} \\
                          \end{array}
                        \right),\; with\; D(\Lambda_{1}^{\alpha})=X^{\alpha}\times X^{\alpha}\times X^{\alpha}.
                         \end{align*}
\end{theorem}
\proof Let $\alpha\in(0,1)$. From \eqref{p1} and \eqref{p3} of the proposition \eqref{prop1} and \eqref{puiss-alpha}, we have:
\begin{align*}\Lambda_{1}^{-\alpha}=\left(
                           \begin{array}{ccc}
 \frac{\sin(\pi\alpha)}{\pi}\int_{0}^{\infty}s^{-\alpha}(s+A)^{-1}ds & 0 & 0 \\
 0 & \frac{\sin(\pi\alpha)}{\pi}\int_{0}^{\infty}s^{-\alpha}(s+A)^{-1}ds & 0 \\
 -\frac{\sin(\pi\alpha)}{\pi}\int_{0}^{\infty}s^{-\alpha}(s+A)^{-2}ds & 0  &
 \frac{\sin(\pi\alpha)}{\pi}\int_{0}^{\infty}s^{-\alpha}(s+A)^{-1}ds\\
 \end{array}
  \right).\end{align*}
  We have \begin{equation*} \frac{\sin(\pi\alpha)}{\pi}\int_{0}^{\infty}s^{-\alpha}(s+A)^{-1}ds=A^{-\alpha}.\end{equation*}
 Applying Lemma \eqref{lemme1}, we obtain:
 \begin{equation*} -\frac{\sin(\pi\alpha)}{\pi}\int_{0}^{\infty}s^{-\alpha}(s+A)^{-2}=-\alpha A^{-1-\alpha}.\end{equation*}
 This gives \begin{align*}
 \Lambda_{1}^{-\alpha}=\left(
                          \begin{array}{ccc}
                            A^{-\alpha} & 0 & 0 \\
                            0 & A^{-\alpha} & 0 \\
                            -\alpha A^{-\alpha-1} & 0 & A^{-\alpha} \\
                          \end{array}
                        \right),\; with\; \Lambda_{1}^{-\alpha}\in\mathcal{L}(X\times X\times X).
 \end{align*}
  Since $\Lambda_{1}^{-\alpha}$ is injective, therefore invertible, and we have \begin{align*}
 \Lambda_{1}^{\alpha}=\left(
                          \begin{array}{ccc}
                            A^{\alpha} & 0 & 0 \\
                            0 & A^{\alpha} & 0 \\
                            \alpha A^{\alpha-1} & 0 & A^{\alpha} \\
                          \end{array}
                        \right),\; with\; D(\Lambda_{1}^{\alpha})=X^{\alpha}\times X^{\alpha}\times X^{\alpha}.
                         \end{align*}
\qed\\
\begin{Proposition}\label{prop32} Let $A$ a positive operator and \begin{align*}\Lambda_{\left(3, \frac{1}{2}\right)}=\left(\begin{array}{ccc}0 & -I & 0 \\ A & 2 A^{\frac{1}{2}} & 0 \\ 0 & 0 & 2 A^{\frac{1}{2}}\end{array}\right)\end{align*}. Then:\begin{enumerate}
                                                        \item \label{p321}  $\Lambda_{\left(3, \frac{1}{2}\right)}$ is a positive operator and $D\left(\Lambda_{\left(3, \frac{1}{2}\right)}\right)=X^{1} \times X^{\frac{1}{2}} \times X^{\frac{1}{2}},$
                                                        \item \label{p322} $\sigma\left(\Lambda_{\left(3, \frac{1}{2}\right)}\right)=\sigma\left(A^{\frac{1}{2}}\right) \cup \sigma\left(2 A^{\frac{1}{2}}\right),$
                                                        \item \label{p323} $\forall s \in \rho\left(-\Lambda_{2}\right),$
                                                      \end{enumerate}
\begin{align*}
\left(s+\Lambda_{\left(3, \frac{1}{2}\right)}\right)^{-1}=\left(\begin{array}{ccc}
\left(s+2 A^{\frac{1}{2}}\right)\left(s+A^{\frac{1}{2}}\right)^{-2} & \left(s+A^{\frac{1}{2}}\right)^{-2} & 0 \\
-A\left(s+A^{\frac{1}{2}}\right)^{-2} & s\left(s+A^{\frac{1}{2}}\right)^{-2} & 0 \\
0 & 0 & \left(s+2 A^{\frac{1}{2}}\right)^{-1}
\end{array}\right).
\end{align*}
\end{Proposition}
\proof \qed
\begin{theorem} Let $A$ a positive operator and \begin{align*}\Lambda_{\left(3, \frac{1}{2}\right)}=\left(\begin{array}{ccc}0 & -I & 0 \\ A & 2 A^{\frac{1}{2}} & 0 \\ 0 & 0 & 2 A^{\frac{1}{2}}\end{array}\right).\end{align*}
The fractional powers of $\Lambda_{\left(3, \frac{1}{2}\right)}$ can be calculated explicitly. For $\alpha \in(0,1)$, we have:

\begin{align*}
\Lambda_{\left(3, \frac{1}{2}\right)}^{-\alpha}=\left(\begin{array}{ccc}
(1+\alpha) A^{-\frac{\alpha}{2}} & \alpha A^{-\frac{1+\alpha}{2}} & 0 \\
-\alpha A^{\frac{1-\alpha}{2}} & (1-\alpha) A^{-\frac{\alpha}{2}} & 0 \\
0 & 0 & 2^{-\alpha} A^{-\frac{\alpha}{2}}
\end{array}\right)\in \mathcal{L}\left(X^{\frac{1}{2}} \times X \times X\right),
\end{align*}

and

\begin{align*}
\Lambda_{\left(3, \frac{1}{2}\right)}^{\alpha}=\left(\begin{array}{ccc}
(1-\alpha) A^{\frac{\alpha}{2}} & -\alpha A^{\frac{\alpha-1}{2}} & 0 \\
\alpha A^{\frac{1+\alpha}{2}} & (1+\alpha) A^{\frac{\alpha}{2}} & 0 \\
0 & 0 & 2^{\alpha} A^{\frac{\alpha}{2}}
\end{array}\right) \text {, with } D\left(\Lambda_{2}^{\alpha}\right)=X^{\frac{1+\alpha}{2}} \times X^{\frac{\alpha}{2}} \times X^{\frac{\alpha}{2}}.
\end{align*}
\end{theorem}
\proof Let $\alpha \in(0,1)$. From \eqref{p321} and \eqref{p323} of Proposition \eqref{prop32}  and \eqref{puiss-alpha}, we have:

$\Lambda_{\left(3, \frac{1}{2}\right)}^{-\alpha}=$\begin{align*}\left(\begin{array}{ccc}
\frac{\sin(\pi\alpha)}{\pi} \int_{0}^{\infty} s^{-\alpha}(s+2 A^{\frac{1}{2}})(s+A^{\frac{1}{2}})^{-2}ds & \frac{\sin (\pi \alpha)}{\pi}\int_{0}^{\infty}s^{-\alpha}(s+A^{\frac{1}{2}})^{-2}ds& 0 \\
-A\frac{\sin (\pi \alpha)}{\pi} \int_{0}^{\infty}s^{-\alpha}(s+A^{\frac{1}{2}})^{-2}ds & \frac{\sin(\pi \alpha)}{\pi} \int_{0}^{\infty} s^{1-\alpha}(s+A^{\frac{1}{2}})^{-2}ds & 0 \\
0 & 0 & \frac{\sin (\pi\alpha)}{\pi} \int_{0}^{\infty} s^{-\alpha}(s+2 A^{\frac{1}{2}})^{-1}ds
\end{array}\right).
\end{align*}
From Lemma (3.2) we have
\begin{align*}
\frac{\sin (\pi \alpha)}{\pi} \int_{0}^{\infty} s^{-\alpha}(s+2 A^{\frac{1}{2}}(s+A^{\frac{1}{2}})^{-2}ds 
&= \frac{\sin (\pi \alpha)}{\pi}\Big[\int_{0}^{\infty} s^{-\alpha}(s+A^{\frac{1}{2}})^{-1} ds + A^{\frac{1}{2}}\int_{0}^{\infty}s^{-\alpha}(s+A^{\frac{1}{2}})^{-2} d s\Big] \\
&= A^{-\frac{\alpha}{2}} + \alpha A^{\frac{1}{2}} A^{\frac{-1-\alpha}{2}} \\
&= (1+\alpha) A^{-\frac{\alpha}{2}} \\
& \\
\frac{\sin (\pi \alpha)}{\pi} \int_{0}^{\infty} s^{-\alpha}\left(s+A^{\frac{1}{2}}\right)^{-2} d s 
&= \alpha A^{-\frac{1+\alpha}{2}},
\end{align*}
and
\begin{align*}
\frac{\sin (\pi \alpha)}{\pi} \int_{0}^{\infty} s^{1-\alpha}\left(s+A^{\frac{1}{2}}\right)^{-2} d s=(1-\alpha) A^{-\frac{\alpha}{2}}.
\end{align*}
From alternative formula:
\begin{align*}
\frac{\sin (\pi \alpha)}{\pi} \int_{0}^{\infty} s^{-\alpha}\left(s+(4 A)^{\frac{1}{2}}\right)^{-1} d s=(4 A)^{-\frac{\alpha}{2}}=2^{-\alpha} A^{-\frac{\alpha}{2}}.
\end{align*}
This gives:
\begin{align*}
\Lambda_{\left(3, \frac{1}{2}\right)}^{-\alpha}=\left(\begin{array}{ccc}
(1+\alpha) A^{-\frac{\alpha}{2}} & \alpha A^{-\frac{1+\alpha}{2}} & 0 \\
-\alpha A^{\frac{1-\alpha}{2}} & (1-\alpha) A^{-\frac{\alpha}{2}} & 0 \\
0 & 0 & 2^{-\alpha} A^{-\frac{\alpha}{2}}
\end{array}\right), \text { and } \Lambda_{\left(3, \frac{1}{2}\right)}^{-\alpha} \in \mathcal{L}\left(X^{\frac{1}{2}} \times X \times X\right).
\end{align*}

Since $\Lambda_{\left(3, \frac{1}{2}\right)}^{-\alpha}$ is injective, therefore invertible, and we have

\begin{align*}
\Lambda_{\left(3, \frac{1}{2}\right)}^{\alpha}=\left(\begin{array}{ccc}
(1-\alpha) A^{\frac{\alpha}{2}} & -\alpha A^{\frac{\alpha-1}{2}} & 0 \\
\alpha A^{\frac{1+\alpha}{2}} & (1+\alpha) A^{\frac{\alpha}{2}} & 0 \\
0 & 0 & 2^{\alpha} A^{\frac{\alpha}{2}}
\end{array}\right) \text {, with } D\left(\Lambda_{2}^{\alpha}\right)=X^{\frac{1+\alpha}{2}} \times X^{\frac{\alpha}{2}} \times X^{\frac{\alpha}{2}}.
\end{align*}
\qed
\begin{example}
$$
\begin{cases}u_{t}=-v, & x \in \Omega, t>0 \\ v_{t}=-A u+2 A^{\frac{1}{2}} v, & x \in \Omega, t>0 \\ w_{t}=2 A^{\frac{1}{2}} w, & x \in \Omega, t>0, \\ u(x, 0)=u_{0}(x), v(x, 0)=v_{0}(x), w(x, 0)=w_{0}(x), & x \in \Omega, \\ u(x, t)=v(x, t)=w(x, t)=0, & x \in \partial \Omega\end{cases}
$$

Where: $A=-\Delta_{D}$ is the Laplacian with Dirichlet boundary conditions, a positive operator in $L^{2}(\Omega)$.
\end{example}
 \subsection{Using change of variables}
 In some cases: $\widetilde{A}_{ij}(s)=(h(s)+A)^{-1}$, with $A\in\mathcal{P}(X)$ and $h:(0,\infty)\rightarrow \mathbb{R}$, therefore:
 \begin{align} \frac {\sin(\pi\alpha)}{\pi}\int_{0}^{\infty} s^{-\alpha}\widetilde{A}_{ij}(s)ds=\frac {\sin(\pi\alpha)}{\pi}\int_{0}^{\infty} s^{-\alpha}(h(s)+A)^{-1}ds.\end{align}
If $h$ is a diffeomorphism mapping $(0,\infty)$ onto $(0,\infty)$, we can introduce the change of variable $u=h(s)$ and rewrite accordingly
 \begin{align} \frac {\sin(\pi\alpha)}{\pi}\int_{0}^{\infty} s^{-\alpha}\widetilde{A}_{ij}(s)ds=\frac {\sin(\pi\alpha)}{\pi}\int_{\lim\limits_{s \rightarrow 0} h(s)}^{\lim\limits_{s \rightarrow +\infty} h(s)} (h^{- 1}(u))^{-\alpha}(u+A)^{-1}(h^{-1})^{'}(u) du.\end{align} \\  This change of variables allows us, in some cases, to reformulate the right-hand member in terms of
fractional powers of $A$.
\begin{Lemma} If $A\in\mathcal{P}(X)$ and $\widetilde{A}_{ij}(s)=(h(s)+A)^{-1}$ for $s>0$ and $h:(0,\infty)\rightarrow(0,\infty)$ with:
\begin{align*} h(s)=\omega s^{\theta},  \omega>0, \theta>0.\end{align*} Then $h$ is a diffeomorphism and:
\begin{align*} \int_{0}^{\infty} s^{-\gamma}\widetilde{A}_{ij}(s)ds=\int_{0}^{\infty} s^{-\gamma}(h(s)+A)^{-1}(s)ds=\frac{1}{\theta}(\frac{1}{\omega})^{\frac{1-\gamma}{\theta}}\frac{\pi}{\sin(\pi\beta)}A^{-\beta}.\end{align*}
With $\beta=1+\frac{\gamma-1}{\theta}$ et $1-\theta<\gamma<1.$
\label{chvar}
\end{Lemma}
\proof Change of variable $u=\omega s^{\theta}$ gives:
\begin{align*} \int_{0}^{\infty} s^{-\gamma}(h(s)+A)^{-1}(s)ds=\int_{0}^{\infty}(\frac{1}{\omega})^{-\frac{\gamma}{\theta}}u^{-\frac{\gamma}{\theta}}(u+A)^{-1}(\frac{1}{\omega})^{\frac{1}{\theta}}\frac{1}{\theta}u^{\frac{1}{\theta}-1}du
\\=\frac{1}{\theta}(\frac{1}{\omega})^{\frac{1-\gamma}{\theta}}\int_{0}^{\infty}u^{\frac{1-\gamma}{\theta}-1}(u+A)^{-1}du\\
=\frac{1}{\theta}(\frac{1}{\omega})^{\frac{1-\gamma}{\theta}}\frac{\pi}{\sin(\pi\beta)}[\frac{\sin(\pi\beta)}{\pi}\int_{0}^{\infty}u^{-\beta}(u+A)^{-1}du]\end{align*}
With $\beta=1+\frac{\gamma-1}{\theta}$, since $1-\theta<\gamma<1$, then $0<\beta<1$ and the expression in square brackets is exactly $A^{-\beta}$.\qed
\begin{Proposition}\label{lam4} Let $A$ a positive operator and \begin{align*}\Lambda_{4}=\left(\begin{array}{ccc}0 & 0 & -I \\ 0 & A & 0 \\ A & 0 & 0\end{array}\right).\end{align*}
Then:\begin{enumerate}
       \item \label{p41} $\Lambda_{4}$ is a positive operator and $D\left(\Lambda_{4}\right)=X^{1} \times X^{1} \times X^{\frac{1}{2}}$.
       \item \label{p42} $\sigma\left(\Lambda_{4}\right)=\sigma( \pm i A^{\frac{1}{2}}) \cup \sigma(A)$.
       \item \label{p43} $\forall s \in \rho\left(-\Lambda_{4}\right)$ :
       \begin{align*}
\left(s+\Lambda_{4}\right)^{-1}=\left(\begin{array}{ccc}
s\left(s^{2}+A\right)^{-1} & 0 & \left(s^{2}+A\right)^{-1} \\
0 & (s+A)^{-1} & 0 \\
-A\left(s^{2}+A\right)^{-1} & 0 & s\left(s^{2}+A\right)^{-1}
\end{array}\right)
\end{align*}
     \end{enumerate}
\end{Proposition}
\proof \qed
\begin{theorem} let $A$ a positive operator and \begin{align*}\Lambda_{4}=\left(\begin{array}{ccc}0 & 0 & -I \\ 0 & A & 0 \\ A & 0 & 0\end{array}\right)\end{align*}.
Fractional powers of $\Lambda_{4}$ can be calculated explicitly. For $\alpha \in(0,1)$, we have:
\begin{align*}
\Lambda_{4}^{-\alpha}=\left(\begin{array}{ccc}
\cos \left(\frac{\pi \alpha}{2}\right) A^{-\frac{\alpha}{2}} & 0 & \sin \left(\frac{\pi \alpha}{2}\right) A^{-\frac{\alpha+1}{2}} \\
0 & A^{-\alpha} & 0 \\
-\sin \left(\frac{\pi \alpha}{2}\right) A^{\frac{1-\alpha}{2}} & 0 & \cos \left(\frac{\pi \alpha}{2}\right) A^{-\frac{\alpha}{2}}
\end{array}\right) \text {, with } \Lambda_{4}^{-\alpha} \in \mathcal{L}\left(X^{\frac{1}{2}} \times X \times X\right)
\end{align*}
and
\begin{align*}
\Lambda_{4}^{\alpha}=\left(\begin{array}{ccc}
\cos \left(\frac{\pi \alpha}{2}\right) A^{\frac{\alpha}{2}} & 0 & -\sin \left(\frac{\pi \alpha}{2}\right) A^{\frac{\alpha-1}{2}} \\
0 & A^{\alpha} & 0 \\
\sin \left(\frac{\pi \alpha}{2}\right) A^{\frac{1+\alpha}{2}} & 0 & \cos \left(\frac{\pi \alpha}{2}\right) A^{\frac{\alpha}{2}}
\end{array}\right) \text {, with } D\left(\Lambda_{4}^{\alpha}\right)=X^{\frac{1+\alpha}{2}} \times X^{\alpha} \times X^{\frac{\alpha}{2}}.\end{align*}
\end{theorem}

\proof  Let $\alpha \in(0,1)$. From \eqref{p41} and \eqref{p43} of the proposition\eqref{lam4} and \eqref{puiss-alpha}, we have:
\begin{align*}\Lambda_{4}^{-\alpha}=\left(\begin{array}{ccc}\frac{\sin (\pi \alpha)}{\pi} \int_{0}^{\infty} s^{1-\alpha}\left(s^{2}+A\right)^{-1} d s & 0 & \frac{\sin (\pi \alpha)}{\pi} \int_{0}^{\infty}s^{-\alpha}\left(s^{2}+A\right)^{-1} d s \\ 0 & \frac{\sin (\pi \alpha)}{\pi} \int_{0}^{\infty} s^{-\alpha}(s+A)^{-1} d s & 0 \\-A \frac{\sin (\pi \alpha)}{\pi} \int_{0}^{\infty} s^{-\alpha}\left(s^{2}+A\right)^{-1} d s&0&\frac{\sin (\pi \alpha)}{\pi} \int_{0}^{\infty} s^{1-\alpha}\left(s^{2}+A\right)^{-1} d s\end{array}\right).\end{align*}
We have
\begin{align*}
\frac{\sin (\pi \alpha)}{\pi} \int_{0}^{\infty} s^{-\alpha}(s+A)^{-1} d s=A^{-\alpha}.
\end{align*}

Using Lemma \eqref{chvar} with $\theta=2, \gamma=\alpha-1, \beta=\frac{1+\gamma}{2}=\frac{\alpha}{2}$ and $\omega=1$ :
\begin{align*}
\frac{\sin (\pi \alpha)}{\pi} \int_{0}^{\infty} s^{1-\alpha}\left(s^{2}+A\right)^{-1} d s=\cos \left(\frac{\alpha \pi}{2}\right) A^{\frac{-\alpha}{2}}.
\end{align*}

Using Lemma(3.8) with $\theta=2, \gamma=\alpha, \beta=\frac{1+\gamma}{2}=\frac{\alpha+1}{2}$ and $\omega=1$ :
\begin{align*}
\frac{\sin (\pi \alpha)}{\pi} \int_{0}^{\infty} s^{-\alpha}\left(s^{2}+A\right)^{-1} d s=\sin \left(\frac{\alpha \pi}{2}\right) A^{-\frac{\alpha+1}{2}}
\end{align*}
and
\begin{align*}
-A \frac{\sin (\pi \alpha)}{\pi} \int_{0}^{\infty} s^{-\alpha}\left(s^{2}+A\right)^{-1} d s=-\sin \left(\frac{\alpha \pi}{2}\right) A^{\frac{1-\alpha}{2}}.
\end{align*}

This gives:
\begin{align*}
\Lambda_{4}^{-\alpha}=\left(\begin{array}{ccc}
\cos \left(\frac{\pi \alpha}{2}\right) A^{-\frac{\alpha}{2}} & 0 & \sin \left(\frac{\pi \alpha}{2}\right) A^{-\frac{\alpha+1}{2}} \\
0 & A^{-\alpha} & 0 \\
-\sin \left(\frac{\pi \alpha}{2}\right) A^{\frac{1-\alpha}{2}} & 0 & \cos \left(\frac{\pi \alpha}{2}\right) A^{-\frac{\alpha}{2}}
\end{array}\right) \text {, with } \Lambda_{4}^{-\alpha} \in \mathcal{L}\left(X^{\frac{1}{2}} \times X \times X\right)
\end{align*}

Since $\Lambda_{4}^{-\alpha}$ is injective, therefore invertible, and we have
\begin{align*}
\Lambda_{4}^{\alpha}=\left(\begin{array}{ccc}
\cos \left(\frac{\pi \alpha}{2}\right) A^{\frac{\alpha}{2}} & 0 & -\sin \left(\frac{\pi \alpha}{2}\right) A^{\frac{\alpha-1}{2}} \\
0 & A^{\alpha} & 0 \\
\sin \left(\frac{\pi \alpha}{2}\right) A^{\frac{1+\alpha}{2}} & 0 & \cos \left(\frac{\pi \alpha}{2}\right) A^{\frac{\alpha}{2}}
\end{array}\right) \text {, with } D\left(\Lambda_{4}^{\alpha}\right)=X^{\frac{1+\alpha}{2}} \times X^{\alpha} \times X^{\frac{\alpha}{2}}
\end{align*}\qed
\begin{example} Consider the following coupled reaction-diffusion System with Dirichlet Boundary Conditions:
\begin{align*}
\begin{cases}u_{t}=-w, & x \in \Omega, t>0  \tag{16}\\ v_{t}=-\Delta v, & x \in \Omega, t>0 \\ w_{t}=-\Delta u, & x \in \Omega, t>0 \\ u(x, t)=v(x, t)=w(x, t)=0, & x \in \partial \Omega, t>0 \\ u(x, 0)=u_{0}(x), v(x, 0)=v_{0}(x), w(x, 0)=w_{0}(x), & x \in \Omega\end{cases}
\end{align*}

Note that the linear operator associated to this system has the same format of $\Lambda_{4}$ given in \eqref{lam4}.
\end{example}
The following figure illustrates the location of the spectra of $\Lambda_{4}$ and its fractional power for $\alpha \in (0,1)$, taking $A = -\Delta_{D}$, which is a positive operator with a discrete spectrum represented by an increasing sequence of strictly positive real numbers $(\mu_{n})_{n\in\mathbb{N}^{*}}$.
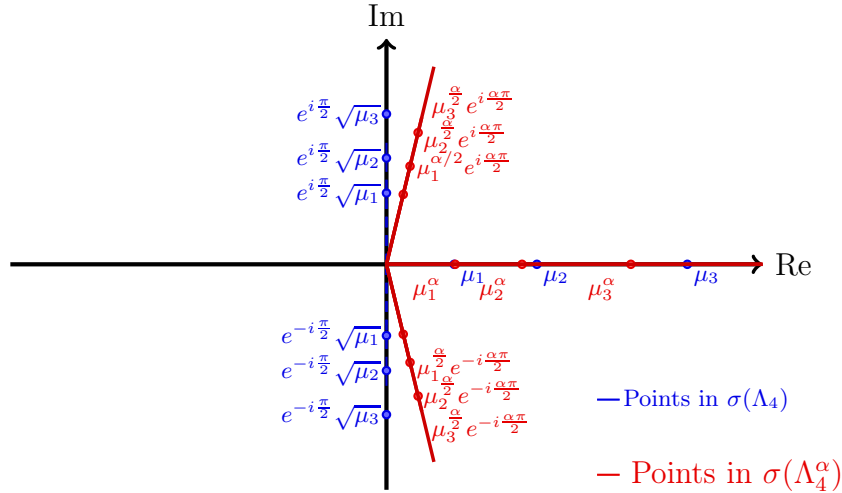
\begin{figure}[ht]
\centering
\begin{tikzpicture}[scale=1]
    \draw[ultra thick, ->] (-5,0) -- (5,0) node[right]{$\mathrm{Re}$};
    \draw[ultra thick, ->] (0,-3) -- (0,3) node[above]{$\mathrm{Im}$};
    
    \def\muOne{0.9} \def\muTwo{2} \def\muThree{4}
    \def\alphaVal{0.85}
    \def\pointSize{1.4pt}

    \tikzset{
        point bleu/.style={blue!90!black, thick, fill=blue!60!white},
        label bleu/.style={blue!90!black, font=\scriptsize, inner sep=1pt},
        point rouge/.style={red!90!black, thick, fill=red!60!white},
        label rouge/.style={red!90!black},
        ligne bleue/.style={blue!80!black, thick, dash pattern=on 4pt off 3pt},
        ligne rouge/.style={red!80!black, thick, line width=1.3pt}
    }

    \foreach \n/\x in {1/\muOne, 2/\muTwo, 3/\muThree} {
        \filldraw[point bleu] (\x,0) circle (\pointSize)
            node[label bleu, below right=1pt] {$\mu_{\n}$};
        
        \filldraw[point bleu] (0,{sqrt(\x)}) circle (\pointSize)
            node[label bleu, left=1pt] {$e^{i\frac{\pi}{2}}\sqrt{\mu_{\n}}$};
        \filldraw[point bleu] (0,{-sqrt(\x)}) circle (\pointSize)
            node[label bleu, left=1pt] {$e^{-i\frac{\pi}{2}}\sqrt{\mu_{\n}}$};
        
        \draw[ligne bleue] (0,0) -- (0.85*\x,0);
        \draw[ligne bleue] (0,0) -- (0,{0.85*sqrt(\x)});
        \draw[ligne bleue] (0,0) -- (0,{-0.85*sqrt(\x)});
    }

    \foreach \n/\x in {1/\muOne, 2/\muTwo, 3/\muThree} {
        \pgfmathsetmacro{\radA}{pow(\x,\alphaVal)}
        \pgfmathsetmacro{\radB}{pow(\x,\alphaVal/2)}
        \pgfmathsetmacro{\angle}{deg(\alphaVal*pi/2)}
        
        \filldraw[point rouge] (\radA,0) circle (\pointSize)
            node[label rouge, below left=1pt, font=\scriptsize] {$\mu_{\n}^{\alpha}$};
        
        \filldraw[point rouge] ({\radB*cos(\angle)}, {\radB*sin(\angle)}) circle (\pointSize);
        \filldraw[point rouge] ({\radB*cos(-\angle)}, {\radB*sin(-\angle)}) circle (\pointSize);
        
        \ifthenelse{\n=1}{
            \node[label rouge, above right=1pt, font=\tiny] at ({\radB*cos(\angle)}, {\radB*sin(\angle)}) 
                {$\mu_{\n}^{\alpha/2}e^{i\frac{\alpha\pi}{2}}$};
            \node[label rouge, below right=1pt, font=\tiny] at ({\radB*cos(-\angle)}, {\radB*sin(-\angle)}) 
                {$\mu_{\n}^{\frac{\alpha}{2}}e^{-i\frac{\alpha\pi}{2}}$};
        }{
            \node[label rouge, above right=1pt, font=\scriptsize] at ({\radB*cos(\angle)}, {\radB*sin(\angle)}) 
                {$\mu_{\n}^{\frac{\alpha}{2}}e^{i\frac{\alpha\pi}{2}}$};
            \node[label rouge, below right=1pt, font=\scriptsize] at ({\radB*cos(-\angle)}, {\radB*sin(-\angle)}) 
                {$\mu_{\n}^{\frac{\alpha}{2}}e^{-i\frac{\alpha\pi}{2}}$};
        }
        
        \draw[ligne rouge] (0,0) --+ (5,0); 
        \draw[ligne rouge] (0,0) --+ ({\radB*cos(\angle)*1.5}, {\radB*sin(\angle)*1.5}); 
        \draw[ligne rouge] (0,0) --+ ({\radB*cos(-\angle)*1.5}, {\radB*sin(-\angle)*1.5}); 
    }

    \begin{scope}[shift={(2.8,-1.8)}]
        \draw[point bleu] (0,0) -- +(0.3,0) node[label bleu, right]{Points in $\sigma(\Lambda_{4})$};
        \draw[point rouge] (0,-1) -- +(0.3,0) node[label rouge, right]{Points in $\sigma(\Lambda_{4}^{\alpha})$};
    \end{scope}
\end{tikzpicture}
\caption{\textbf{Spectrum of $\Lambda_4$ and $\Lambda_4^\alpha$ in the complex plane}}
\label{figure_spectre}
\end{figure}

\subsection{Using the second resolvent identity}
  \noindent There are situations where the $\widetilde{A}_{ij}(s)$ in the expression \eqref{hyp} are of the form:\\
\begin{align*}\widetilde{A}_{ij}(s)=(s+A_{1})^{-1}(s+A_{2})^{-1} \end{align*}
\begin{Lemma}
Let $A_{1}$ and $A_{2}$ positives operators in $X$.\\ If $A_{1}-A_{2}:D(A_{1})\cap D(A_{2})\subset X\rightarrow X$ is invertible and commutes with the resolvent of $A_{2}$, then :\begin{enumerate}
                                            \item $\forall\lambda\in\rho(-A_{1})\cap\rho(-A_{2})$,
                                            \begin{equation*}(\lambda+A_{2})^{-1}(\lambda+A_{1})^{-1}=(A_{1}-A_{2})^{-1}[(\lambda+A_{2})^{-1}-(\lambda+A_{1})^{-1}].
                                            \end{equation*}
                                            \item \begin{equation*}\frac{\sin(\pi\alpha)}{\pi}\int_{0}^{\infty}s^{-\alpha}(s+A_{2})^{-1}(s+A_{1})^{-1}ds=(A_{1}-A_{2})^{-1}[A_{2}^{-\alpha}-A_{1}^{-\alpha}].
                                                \end{equation*}\\
                                         \end{enumerate}
                                          \label{lemmeresolv}
\end{Lemma}
\proof
$\forall\lambda\in\rho(-A_{1})\cap\rho(-A_{2})$:
\begin{align*}
(\lambda+A_{2})^{-1}-(\lambda+A_{1})^{-1}
  &=(\lambda+A_{2})^{-1}(\lambda+A_{1})(\lambda+A_{1})^{-1}-(\lambda+A_{2})^{-1}(\lambda+A_{2})(\lambda+A_{1})^{-1}\\
  &=(\lambda+A_{2})^{-1}(\lambda+A_{1}-\lambda-A_{2})(\lambda+A_{1})^{-1}\\
  &=(\lambda+A_{2})^{-1}(A_{1}-A_{2})(\lambda+A_{1})^{-1}.
  \end{align*}
  Since the resolvent of $A_{2}$ commutes with $A_{1}-A_{2}$ which is invertible, then:\\
  $(\lambda+A_{2})^{-1}(\lambda+A_{1})^{-1}=(A_{1}-A_{2})^{-1}[(\lambda+A_{2})^{-1}-(\lambda+A_{1})^{-1}]$.
  Since $[0,\infty)\subset\rho(-A_{1})\cap\rho(-A_{2})$, we can use the previous equation to obtain:\\
 \begin{align*}
 &\frac{\sin(\pi\alpha)}{\pi}\int_{0}^{\infty}s^{-\alpha}(s+A_{2})^{-1}(s+A_{1})^{-1}ds\\
 &=(A_{1}-A_{2})^{-1}\frac{\sin(\pi\alpha)}{\pi}\int_{0}^{\infty}(s^{-\alpha}(s+A_{2})^{-1}-s^{-\alpha}(s+A_{1})^{-1})ds\\               &=(A_{1}-A_{2})^{-1}[A_{2}^{-\alpha}-A_{1}^{-\alpha}].
                                                \end{align*}
\qed
\begin{example} Weakly coupled reaction-diffusion system:\\
We consider the following DPE, with $a_{1}$, $a_{2}$ and $a_{3}$ being strictly positive and distinct in pairs:
\begin{equation}
\begin{cases}
      u_{t}-a_{1}\Delta_{D}u=f,   x\in\Omega, t>0,\\
      v_{t}-a_{2}\Delta_{D}v=g,  x\in\Omega, t>0,\\
      w_{t}+u-a_{3}\Delta_{D}w=g,  x\in\Omega, t>0,\\
      u(x,0)=v(x,0)=0=w(x,0),             x\in\partial\Omega.
    \end{cases}\,.
\label{edp3}
\end{equation}
 If $A_{i}=-a_{i}\Delta_{D}$ (positive operators), then system \eqref{edp3} is written:
 \begin{align}
 \frac{d}{dt}\left(
               \begin{array}{c}
                 u \\
                 v \\
                 w \\
               \end{array}
             \right)+\left(
                       \begin{array}{ccc}
                         A_{1} & 0 & 0 \\
                         0 & A_{2} & 0 \\
                         I & 0 & A_{3} \\
                       \end{array}
                     \right)\left(
                              \begin{array}{c}
                                u \\
                                v \\
                                w \\
                              \end{array}
                            \right)=\left(
                                      \begin{array}{c}
                                        f \\
                                        g \\
                                        h \\
                                      \end{array}
                                    \right).
 \end{align}
\end{example}
\begin{Lemma} Let $A_{1}$, $A_{2}$ , $A_{3}$ positive operators with same domain $X^{1}$ and $\Lambda_{3}=\left(
                       \begin{array}{ccc}
                         A_{1} & 0 & 0 \\
                         0 & A_{2} & 0 \\
                         I & 0 & A_{3} \\
                       \end{array}
                     \right)$. \begin{align*}(\Lambda_{3},X^{1}\times X^{1}\times X^{1}) \mbox{ is closed with dense domain}.\end{align*}\end{Lemma}
   \proof Since $\overline{X^{1}}=X$, then $\overline{D\left(\Lambda_{3}\right)}=X \times X \times X$. Let

\begin{align*}
U_{n}=\left(\begin{array}{c}
u_{n} \\
v_{n} \\
w_{n}
\end{array}\right) \in D\left(\Lambda_{3}\right)^{\mathbb{N}}
\end{align*}

such that $U_{n} \rightarrow U=\left(\begin{array}{c}u \\ v \\ w\end{array}\right)$ and $\Lambda_{3} U_{n} \rightarrow V=\left(\begin{array}{l}x \\ y \\ z\end{array}\right)$. We have

\begin{align*}
\left\{\begin{array}{l}
u_{n} \rightarrow u \quad \text { and } \quad A_{1} u_{n} \rightarrow x \\
v_{n} \rightarrow v \quad \text { and } A_{2} v_{n} \rightarrow y \\
w_{n} \rightarrow w \quad \text { and } A_{3} w_{n} \rightarrow z-x
\end{array}\right.
\end{align*}
Since $A_{1}, A_{2}, A_{3}$ are closed, then $u, v, w \in X^{1}$ and

\begin{align*}
\left\{\begin{array}{l}
A_{1} u=x \\
A_{2} v=y \\
x-A_{3} w=z
\end{array}\right.
\end{align*}

Thus, $U \in D\left(\Lambda_{3}\right)$ and $\Lambda_{3} U=V$. Therefore, $\Lambda_{3}$ is closed.\qed \\The results of the following proposition can be easly proved. 
   \begin{Proposition}
   Let $A_{1}$, $A_{2}$ , $A_{3}$ positive operators with same domain $X^{1}$ and \begin{align*}\Lambda_{3}=\left(
                       \begin{array}{ccc}
                         A_{1} & 0 & 0 \\
                         0 & A_{2} & 0 \\
                         I & 0 & A_{3} \\
                       \end{array}
                     \right).\end{align*} If $A_{i}-A_{j}$ is invertible and commutes with the resolvent of $A_{k}$ for any pairwise $i$, $j$ and $k$, then
 \begin{enumerate}
     \item \label{f1} $\Lambda_{3}$ is positive,
     \item \label{f2}$\sigma(\Lambda_{3})=\sigma(A_{1})\cup\ \sigma(A_{2})\cup\sigma(A_{3})$,
     \item \label{f3}\begin{align*}
     (\lambda+\Lambda_{3})^{-1}=\left(
                              \begin{array}{ccc}
                                (\lambda+A_{1})^{-1} & 0 & 0 \\
                                0 & (\lambda+A_{2})^{-1} & 0 \\
                                -(\lambda+A_{1})^{-1}(\lambda+A_{3})^{-1} & 0 & (\lambda+A_{3})^{-1} \\
                              \end{array}
                            \right),\; \forall\lambda\in\rho(-\Lambda_{3}).
                            \end{align*}
   \end{enumerate}\label{finale}
   \end{Proposition}
 \begin{theorem}
  Fractional powers of $\Lambda_{3}$ can be calculated explicitly.\\
  for all $\alpha\in(0,1)$, we have:
 \begin{equation*}\Lambda_{3}^{-\alpha}=\left(
                         \begin{array}{ccc}
                           A_{1}^{-\alpha} & 0 & 0 \\
                           0 & A_{2}^{-\alpha} & 0 \\
                           (A_{1}-A_{3})^{-1}(A_{1}^{-\alpha}-A_{3}^{-\alpha}) & 0 & A_{3}^{-\alpha} \\
                         \end{array}
                       \right)\in\mathcal{L}(X\times X\times X),
\end{equation*}\\and\\
\begin{equation*}\Lambda_{3}^{\alpha}=\left(
                         \begin{array}{ccc}
                           A_{1}^{\alpha} & 0 & 0 \\
                           0 & A_{2}^{\alpha} & 0 \\
                           (A_{1}-A_{3})^{-1}(A_{1}^{\alpha}-A_{3}^{\alpha}) & 0 & A_{3}^{\alpha} \\
                         \end{array}
                       \right),\;with\; D(\Lambda_{3}^{\alpha})=X^{\alpha}\times X^{\alpha}\times X^{\alpha}.
\end{equation*}
\end{theorem}
\proof Let $\alpha\in(0,1)$. Using \eqref{f3} from proposition \eqref{finale} and lemma \eqref{lemmeresolv}, we obtain \begin{equation*}\Lambda_{3}^{-\alpha}=\left(
                         \begin{array}{ccc}
                           A_{1}^{-\alpha} & 0 & 0 \\
                           0 & A_{2}^{-\alpha} & 0 \\
                           (A_{1}-A_{3})^{-1}(A_{1}^{-\alpha}-A_{3}^{-\alpha}) & 0 & A_{3}^{-\alpha} \\
                         \end{array}
                       \right)\in\mathcal{L}(X\times X\times X).
\end{equation*} $\Lambda_{3}^{-\alpha}$ is injective, thus invertible and \begin{equation*}\Lambda_{3}^{\alpha}=\left(
                         \begin{array}{ccc}
                           A_{1}^{\alpha} & 0 & 0 \\
                           0 & A_{2}^{\alpha} & 0 \\
                           (A_{1}-A_{3})^{-1}(A_{1}^{\alpha}-A_{3}^{\alpha}) & 0 & A_{3}^{\alpha} \\
                         \end{array}
                       \right),\;with\; D(\Lambda_{3}^{\alpha})=X^{\alpha}\times X^{\alpha}\times X^{\alpha}.
\end{equation*} \qed
\section*{Acknowledgments}  

The authors would like to thank Nassim Athmouni for his valuable insights and constructive discussions.  

\end{document}